\documentclass[12pt,thmsa]{article}
\usepackage{amssymb}
\usepackage{graphicx}

\begin{document}

\author{Francis OGER}
\title{Elementary equivalence\\
for abelian-by-finite and nilpotent groups}
\date{}

\begin{center}
\textbf{Self-avoiding space-filling}

\textbf{folding curves in dimension }$3$

Francis OGER\bigskip
\end{center}

\bigskip

\noindent \textbf{Abstract. }Various examples of folding curves in $%
\mathbb{R}
^{2}$ have been considered: dragons and other square curves, terdragons and
other triangular curves, Peano-Gosper curves based on hexagons. They are
self-avoiding. They form coverings of $%
\mathbb{R}
^{2}$, by one curve or by a small number of curves, which satisfy the local
isomorphism property. They were used to define some fractals. We construct
an example with similar properties in $%
\mathbb{R}
^{3}$.\bigskip

Various examples of folding curves in $%
\mathbb{R}
^{2}$ have been considered:

\noindent - square folding curves such as the dragon curve are constructed
from a regular tiling of $%
\mathbb{R}
^{2}$ by squares: $s$-folding curves are obtained by folding $s$ times a
strip of paper in $2$, each time possibly left or right, then unfolding it
with $\pi /2$ angles;

\noindent - triangular folding curves such as the terdragon curve are
constructed from a regular tiling of $%
\mathbb{R}
^{2}$ by equilateral triangles: $s$-folding curves are obtained by folding $%
s $ times a strip of paper in $3$, each time possibly left then right or
right then left, then unfolding it with $\pi /3$ angles;

\noindent - Peano-Gosper curves are constructed from a regular tiling of $%
\mathbb{R}
^{2}$ by hexagons.

These curves have two characteristic properties: First, they are
self-avoiding. Second, $s$-curves of the same type completely cover larger
and larger parts of $%
\mathbb{R}
^{2}$ when $s$ becomes larger. They were used to construct fractals.

Any inductive limit of $s$-curves is a curve with $1$ endpoint or without
endpoint. In [6], [7] and [8], we proved that any such curve can be
completed, in an essentially unique way, into a covering of $%
\mathbb{R}
^{2}$ by a set of curves without endpoints which satisfies the local
isomorphism property defined below.

This covering consists of at most $6$ curves. It is locally isomorphic to a
covering constructed in the same way which consists of only $1$ curve.

The object of the present paper is to construct a similar example in $%
\mathbb{R}
^{3}$. As only regular tilings by cubes exist in $%
\mathbb{R}
^{3}$, it is natural to use them. For that purpose, one idea was to fold
several times a wire in $2$ according to the $3$ possible directions of the
space, then unfold it with $\pi /2$ angles. However, [4] shows that the
curves obtained in that way are generally not self-avoiding. Another
possible direction of research was mentioned in [1].

In the first section, we show how convolution products can be used to
construct coverings of $%
\mathbb{R}
^{2}$, and also coverings of $%
\mathbb{R}
^{n}$\ for $n\geq 3$, by sets of folding curves which satisfy the local
isomorphism property. In the second section, we give an example in $%
\mathbb{R}
^{3}$.\bigskip

\textbf{1. Coverings constructed from convolution products.}\bigskip

For each integer $n\geq 2$, we consider $%
\mathbb{R}
^{n}$\ equipped with the euclidian distance $\delta $. For each $x\in 
\mathbb{R}
^{n}$ and each $s\in 
\mathbb{R}
^{+}$, we write $B(x,s)=\{y\in 
\mathbb{R}
^{n}\mid \delta (x,y)\leq s\}$.

For any $E,F\subset 
\mathbb{R}
^{n}$, an \emph{isomorphism} from $E$ to $F$ is a translation $\tau $ such
that $\tau (E)=F$.\ They are \emph{locally isomorphic} if, for each $x\in 
\mathbb{R}
^{n}$\ (resp. $y\in 
\mathbb{R}
^{n}$) and each $s\in 
\mathbb{R}
^{+}$, there exists $y\in 
\mathbb{R}
^{n}$\ (resp. $x\in 
\mathbb{R}
^{n}$) such that $(B(x,s)\cap E,x)\cong (B(y,s)\cap F,y)$.

We say that $E\subset 
\mathbb{R}
^{n}$ satisfies the \emph{local isomorphism property} if, for any $x\in 
\mathbb{R}
^{n}$ and $s\in 
\mathbb{R}
^{+}$, there exists $t\in 
\mathbb{R}
^{+}$ such that each $B(y,t)$ with $y\in 
\mathbb{R}
^{n}$\ contains some $z$ with $(B(z,s)\cap E,z)\cong (B(x,s)\cap E,x)$.

These notions are\ used for the study of aperiodic tilings. There exist some
finite sets of prototiles and matching rules from which it is possible to
construct non periodic tilings, but not periodic tilings. Such tilings are
said to be aperiodic.

The first examples in $%
\mathbb{R}
^{2}$ used a large number of prototiles. Later on, R. Penrose constructed an
example with $2$ prototiles. More recently, an example with $1$ prototile
was obtained (see [3]). Some examples were also constructed in $%
\mathbb{R}
^{n}$\ for $n\geq 3$.

In each example, $2^{\omega }$ isomorphism classes of locally isomorphic
tilings are obtained and each tiling satisfies the local isomorphism
property.

Let us consider the tilings of $%
\mathbb{R}
^{n}$\ which are obtained from some given set of prototiles and matching
rules, and the langage which consists of one relational symbol for each
possible configuration of one tile and its neighbours in such a tiling.
Then, by [5, Th. 2.3], two such tilings are locally isomorphic if and only
if the associated relational structures are elementarily equivalent.

The following method is used to construct a non periodic tiling in such a
situation: Each prototile is cut into smaller tiles which are obtained by
reducing the sizes of the prototiles by the same ratio. In that way, from
each bounded cluster of nonoverlapping tiles, we obtain a cluster with more
tiles which satisfies the same rules. By repeting the process, we cover the
whole space.

This operation is called \emph{deflation} and the inverse op\'{e}ration is
called \emph{inflation}. We are going to use similar operations for folding
curves.

We consider the canonical basis $(e_{1},\ldots ,e_{n})$ of $\mathbb{R}^{n}$.
We call \emph{segments} the non oriented intervals $\left[ x,x+e_{i}\right] $
for $x\in 
\mathbb{Z}
^{n}$\ and $1\leq i\leq n$.

A \emph{curve} is a nonoriented sequence of segments which is finite, or
indexed by $%
\mathbb{N}
$, or indexed by $%
\mathbb{Z}
$, and such that:

\noindent (i) two consecutive segments have one and only one common endpoint;

\noindent (ii) the common endpoints of a segment with the preceding one and
with the following one, when these two segments exist, are different.

Two sets of curves are said to be \emph{equivalent} if they are equivalent
modulo a positive isometry, and \emph{isomorphic} if they are equivalent up
to translation.

For any sets of curves $\mathcal{E},\mathcal{F}$, we write $\mathcal{E}\prec 
\mathcal{F}$\ if each segment of a curve of $\mathcal{E}$ is a segment of a
curve of $\mathcal{F}$ and if each connexion between $2$ segments of a curve
of $\mathcal{E}$ also exists in $\mathcal{F}$.

We say that a curve is \emph{self-avoiding} if each of its segments appears
only once. Two curves are said to be \emph{disjoint} if they have no common
segment.

Two consecutive segments of a curve are aligned or form a $\pi /2$ angle. In 
$\mathbb{R}^{2}$, if a curve without aligned consecutive segments is
self-avoiding according to the definition above, then we can make it
self-avoiding in the geometrical sense by rounding the angles between
consecutive segments. For $n\geq 3$, we can do the same thing in $\mathbb{R}%
^{n}$, even for curves with aligned consecutive segments.

In order to preserve the properties of invariance through translations, we
have to do the rounding in the same way for all endpoints of segments such
that the connexions between the segments at\ these endpoints are the same up
to translation.

For $n\geq 2$, we say that a set $\mathcal{E}$ of self-avoiding curves \emph{%
covers} $X\subset 
\mathbb{R}
^{n}$ if each segment containing a point of $X$ appears in one and only one
of the curves. We say that $\mathcal{E}$ is a \emph{covering} of $%
\mathbb{R}
^{n}$ if it covers $%
\mathbb{R}
^{n}$.

A covering of $%
\mathbb{R}
^{n}$ can be interpreted as a tiling by tiles equipped with drawings: we
associate to each segment its Voronoi tile, which consists of the points for
which this segment is the nearest one.

For $k\in 
\mathbb{N}
^{\ast }$ and $l\in k%
\mathbb{N}
^{\ast }$, we say that a set $\mathcal{A}$ of curves of $\mathbb{R}^{n}$
satisfies (H) for $\left( k,l\right) $ if $\mathcal{A}$ is a covering of $%
\mathbb{R}^{n}$, if each cuve of $\mathcal{A}$ has endpoints $x,x+ke_{i}$\
with $x\in k%
\mathbb{Z}
^{n}$\ and $1\leq i\leq n$, if any $2$ such points are joined by a unique
curve of $\mathcal{A}$ and if $l%
\mathbb{Z}
^{n}+\mathcal{A}=\mathcal{A}$.

For each pair $\left( k,l\right) $, each set $\mathcal{A}$ of curves which
satisfies (H) for $\left( k,l\right) $\ and each set $\mathcal{B}$ of
disjoint self-avoiding curves, the\emph{\ convolution} \emph{product} $%
\mathcal{A}\ast \mathcal{B}$ is defined as follows: We apply to $\mathcal{B}$
a homothety of center $0$ and ratio $k$, then we remplace each of the
segments obtained in that way by the curve of $\mathcal{A}$ which has the
same endpoints.

We have $\mathcal{A\prec A}\ast \mathcal{B}$ and $\mathcal{A}\ast \mathcal{B}
$ is a set of disjoint self-avoiding curves. If $\mathcal{B}$ is a covering
(resp. consists of $1$ curve), then $\mathcal{A}\ast \mathcal{B}$ is a
covering (resp. consists of $1$ curve).

If $\mathcal{A}$ and $\mathcal{B}$ satisfy (H) for $\left( i,j\right) $ and $%
\left( k,l\right) $ respectively, then $\mathcal{A}\ast \mathcal{B}$
satisfies (H) for $\left( ik,\mathrm{lcm}(j,il)\right) $, and therefore
satisfies (H) for $\left( ik,jl\right) $; we have $(\mathcal{A}\ast \mathcal{%
B})\ast \mathcal{C}=\mathcal{A}\ast (\mathcal{B}\ast \mathcal{C})$ for each
set $\mathcal{C}$ of disjoint self-avoiding curves.

This product is analogous to the convolution product defined by M. Dekking
in [2, p. 21] for folding sequences in $%
\mathbb{R}
^{2}$. Its two-sided identity is the covering which consists of all curves
with $1$ segment.

For each set $\mathcal{A}$ of curves which satisfies (H) for a pair $\left(
k,l\right) $, we denote by $\mathcal{A}^{0}$ the identity and we write $%
\mathcal{A}^{s+1}=\mathcal{A}\ast \mathcal{A}^{s}$ for each $s\in 
\mathbb{N}
$. Each $\mathcal{A}^{s}$ satisfies (H) for $\left( k^{s},l^{s}\right) $. We
have $\mathcal{A}^{s+t}=\mathcal{A}^{s}\ast \mathcal{A}^{t}$ for any $s,t\in 
\mathbb{N}
$ since $\ast $ is associative.

For each $s\in 
\mathbb{N}
$, we have $\mathcal{A}^{s}\prec \mathcal{A}^{s+1}=\mathcal{A}^{s}\ast 
\mathcal{A}$. We denote by $\widehat{\mathcal{A}}$ the inductive limit of
the sets $\mathcal{A}^{s}$ for $\prec $. In $\widehat{\mathcal{A}}$, the
connexions between segments are defined in each point except $0$. The
equalities $\mathcal{A}^{s+1}=\mathcal{A}\ast \mathcal{A}^{s}$ imply $%
\widehat{\mathcal{A}}=\mathcal{A}\ast \widehat{\mathcal{A}}$.

We write $\overline{e}_{i}=-e_{i}$ for $1\leq i\leq n$\ and\ $E=\{e_{1},%
\overline{e}_{1},\ldots ,e_{n},\overline{e}_{n}\}$. For each pairing
relation $R$ defined on $E$, we denote by $\widehat{\mathcal{A}}_{R}$ the
covering of $\mathbb{R}^{n}$ by curves without endpoints which is obtained
from $\widehat{\mathcal{A}}$ by connecting the segments with endpoint $0$
according to $R$.

For each $x=$ $(x_{1},\ldots ,x_{n})\in 
\mathbb{Z}
^{n}$ and each $h\in 
\mathbb{N}
^{\ast }$, we write

\noindent $\mathrm{P}_{x,h}=[x_{1}-h,x_{1}+h]\times \cdots \times \lbrack
x_{n}-h,x_{n}+h]$.

We say that a covering $\mathcal{E}$ of $\mathbb{R}^{n}$ satisfies (C) if,
for each $x\in 
\mathbb{Z}
^{n}$ and each $h\in 
\mathbb{N}
^{\ast }$, there exists $y\in 
\mathbb{Z}
^{n}-\mathrm{P}_{x,h}$ such that $\mathcal{E}\upharpoonright \mathrm{P}%
_{y,h}\cong \mathcal{E}\upharpoonright \mathrm{P}_{x,h}$. We observe that
the last property is true for each $x\in 
\mathbb{Z}
^{n}$ if it is true for some $w\in 
\mathbb{Z}
^{n}$, since each $\mathrm{P}_{x,h}$ is contained in some $\mathrm{P}_{w,g}$.

For each covering of $\mathbb{R}^{n}$, the local isomorphism property
implies (C). The sets $\widehat{\mathcal{A}}_{R}$ which satisfy (C) are
locally isomorphic since they only differ in $0$.\bigskip

\noindent \textbf{Proposition 1.1}. For each $\widehat{\mathcal{A}}_{R}$,
(C) implies the local isomorphism property.\bigskip

\noindent \textbf{Proof}. For each $h\in 
\mathbb{N}
^{\ast }$, let us consider $x_{h}\in 
\mathbb{Z}
^{n}-\mathrm{P}_{0,h}$ such that $\widehat{\mathcal{A}}_{R}\upharpoonright 
\mathrm{P}_{0,h}\cong \widehat{\mathcal{A}}_{R}\upharpoonright \mathrm{P}%
_{x_{h},h}$. Then there exists $s\in 
\mathbb{N}
^{\ast }$ such that $\mathrm{P}_{x_{h},h}\cap k^{s}%
\mathbb{Z}
^{n}=\emptyset $. For each $y\in x_{h}+l^{s}%
\mathbb{Z}
^{n}$, we have $\mathrm{P}_{y,h}\cap k^{s}%
\mathbb{Z}
^{n}=\emptyset $, and therefore

\noindent $\widehat{\mathcal{A}}_{R}\upharpoonright \mathrm{P}_{y,h}=%
\mathcal{A}^{s}\upharpoonright \mathrm{P}_{y,h}\cong \mathcal{A}%
^{s}\upharpoonright \mathrm{P}_{x_{h},h}=\widehat{\mathcal{A}}%
_{R}\upharpoonright \mathrm{P}_{x_{h},h}\cong \widehat{\mathcal{A}}%
_{R}\upharpoonright \mathrm{P}_{0,h}$.~~$\blacksquare $\bigskip

\noindent \textbf{Remark}. Necessarily, there exists $R$ such that $\widehat{%
\mathcal{A}}_{R}$ satisfies (C) and therefore satisfies the local
isomorphism property. Actually, for each $R$, each $s\in 
\mathbb{N}
^{\ast }$ and each $x\in l^{s}%
\mathbb{Z}
^{n}$, we have $\widehat{\mathcal{A}}_{R}\upharpoonright \mathrm{P}%
_{0,k^{s}-1}\cong \widehat{\mathcal{A}}_{R}\upharpoonright \mathrm{P}%
_{x,k^{s}-1}$ if the connexions of the curves of $\widehat{\mathcal{A}}_{R}$
in $0$ and $x$ are equivalent up to translation.\bigskip

For each set $\mathcal{E}$ of curves, each $x\in 
\mathbb{Z}
^{n}$ and each $s\in 
\mathbb{N}
^{\ast }$, we denote by $\Omega _{s,x}(\mathcal{E})$\ the set of curves
obtained from $\mathcal{E}$ by suppressing the connexions between the
segments which have their common endpoint in $x+k^{s}%
\mathbb{Z}
^{n}$. We have $\Omega _{s,0}(\widehat{\mathcal{A}}_{R})=\mathcal{A}^{s}$
for each $s\in 
\mathbb{N}
^{\ast }$ and each pairing relation $R$ defined on $E$.\bigskip

\noindent \textbf{Proposition 1.2.} No $\widehat{\mathcal{A}}_{R}$ is
invariant through a nontrivial translation.\bigskip

\noindent \textbf{Proof.} We write $\mathcal{E}=\widehat{\mathcal{A}}_{R}$
and we consider $w\in 
\mathbb{Z}
^{n}-\left\{ 0\right\} $ such that $w+\mathcal{E}=\mathcal{E}$.

For each $s\in 
\mathbb{N}
^{\ast }$, we have $l^{s}%
\mathbb{Z}
^{n}+\Omega _{s,0}(\mathcal{E})=\Omega _{s,0}(\mathcal{E})$. If $l^{s}%
\mathbb{Z}
^{n}+\mathcal{E}\neq \mathcal{E}$, then $k^{s}%
\mathbb{Z}
^{n}=\{x\in 
\mathbb{Z}
^{n}\mid l^{s}%
\mathbb{Z}
^{n}+\Omega _{s,x}(\mathcal{E})=\Omega _{s,x}(\mathcal{E})\}$. Consequently, 
$k^{s}%
\mathbb{Z}
^{n}$ is invariant through the translation $v\rightarrow v+w$ and we have $%
w\in k^{s}%
\mathbb{Z}
^{n}$.

For the remainder of the proof, we consider an integer $s$ such that $%
w\notin k^{s}%
\mathbb{Z}
^{n}$, and therefore such that $l^{s}%
\mathbb{Z}
^{n}+\mathcal{E}=\mathcal{E}$. The last property implies $l^{s}%
\mathbb{Z}
^{n}+C\in \mathcal{E}$ for each $C\in \mathcal{E}$. Consequently, the curves
of $\mathcal{E}$ belong to finitely many isomorphism classes.

For each $r\in 
\mathbb{N}
^{\ast }$, each curve $C\in \mathcal{E}$ contains a curve $D\in \mathcal{A}%
^{r}$. For each endpoint $x$ of $D$, there exist some distinct endpoints $%
x_{1},\ldots ,x_{r}$ of segments of $D$ such that $x_{i}\in x+%
\mathbb{Z}
\left\{ e_{1},\ldots ,e_{n}\right\} $\ for $1\leq i\leq r$.

Consequently, for each $C\in \mathcal{E}$ and each $r\in 
\mathbb{N}
^{\ast }$, there exist $i\in \left\{ 1,\ldots ,n\right\} $ and $y_{1},\ldots
,y_{r}$ distinct endpoints of segments of $C$ such that $y_{k}-y_{j}\in 
\mathbb{Z}
e_{i}$\ for $1\leq j,k\leq r$. As each $x\in 
\mathbb{Z}
^{n}$ is an endpoint of only $2n$ segments, it follows that there exist some
distinct segments $S,T$ of $C$ with endpoints $x,y$ such that $y-x\in l^{s}%
\mathbb{Z}
e_{i}$ and $(y-x)+S=T$, which implies $(y-x)+C=C$. Then there exists a
bounded curve $B\subset C$ such that $C=%
\mathbb{Z}
(y-x)+B$.

As the curves $C\in \mathcal{E}$ belong to finitely many isomorphism
classes, we can give a bound for the sizes of the curves $B\subset C$ above.
Using that fact and the properties $\widehat{\mathcal{A}}=\mathcal{A}%
^{r}\ast \widehat{\mathcal{A}}$ for $r\in 
\mathbb{N}
^{\ast }$, we see that the segments of each curve of $\mathcal{E}$ are
aligned, whence a contradiction.~~$\blacksquare $\bigskip

For each set $\mathcal{E}$ of disjoint self-avoiding curves and any $x,y\in 
\mathbb{Z}
^{n}$, we write $x\sim _{\mathcal{E}}y$ if the connexions between segments
of $\mathcal{E}$\ in $x$ and $y$ are the same up to translation.\bigskip

\noindent \textbf{Theorem 1.3}. Let $\mathcal{E}$ be a covering of $\mathbb{R%
}^{n}$ by curves without endpoints which is locally isomorphic to some $%
\widehat{\mathcal{A}}_{R}$. Then there exists a sequence $(x_{s})_{s\in 
\mathbb{N}
^{\ast }}\subset 
\mathbb{Z}
^{n}$, such that $\Omega _{s,x_{s}}(\mathcal{E})\cong \mathcal{A}^{s}$\ for
each $s\in 
\mathbb{N}
^{\ast }$.\bigskip

\noindent \textbf{Proof}. For each $s\in 
\mathbb{N}
^{\ast }$, as elements with the same properties exist for $\widehat{\mathcal{%
A}}_{R}$ by Proposition 1.2, there exist $x,y\in 
\mathbb{Z}
^{n}$ such that $y-x\in l^{s}%
\mathbb{Z}
^{n}$ and $x\nsim _{\mathcal{E}}y$.

For each $h\in 
\mathbb{N}
^{\ast }$ such that $y\in \mathrm{P}_{x,h-1}$, we consider an element $%
z_{h}\in 
\mathbb{Z}
^{n}$ such that the translation $x\rightarrow z_{h}$\ induces an isomorphism 
$\theta _{h}:\mathcal{E}\upharpoonright \mathrm{P}_{x,h}\rightarrow \widehat{%
\mathcal{A}}_{R}\upharpoonright \mathrm{P}_{z_{h},h}$. We have $z_{h}\in
k^{s}%
\mathbb{Z}
^{n}$ since $z_{h}\nsim _{\widehat{\mathcal{A}}_{R}}\theta _{h}(y)$, and $%
\theta _{h}$ induces an isomorphism from $\Omega _{s,x}(\mathcal{E}%
)\upharpoonright \mathrm{P}_{x,h}$\ to $\Omega _{s,z_{h}}(\widehat{\mathcal{A%
}}_{R})\upharpoonright \mathrm{P}_{z_{h},h}=\mathcal{A}^{s}\upharpoonright 
\mathrm{P}_{z_{h},h}$.

There exists an infinite set $M$ of integers such that $y\in \mathrm{P}%
_{x,h-1}$ for $h\in M$ and $z_{g}-z_{h}\in l^{s}%
\mathbb{Z}
^{n}$ for $g,h\in M$. For $g,h\in M$ and $g<h$, the translation $%
z_{g}\rightarrow z_{h}$\ induces an automorphism of $\mathcal{A}^{s}$.
Consequently, each translation $x\rightarrow z_{h}$\ for $h\in M$\ induces
an isomorphism from $\Omega _{s,x}(\mathcal{E})$ to $\mathcal{A}^{s}$. We
write $x_{s}=x$.~~$\blacksquare $\bigskip

\noindent \textbf{Theorem 1.4}. Consider a covering $\mathcal{E}$ of $%
\mathbb{R}^{n}$ by curves without endpoints and a sequence $(x_{s})_{s\in 
\mathbb{N}
^{\ast }}\subset 
\mathbb{Z}
^{n}$ such that $\Omega _{s,x_{s}}(\mathcal{E})\cong \mathcal{A}^{s}$\ for
each $s\in 
\mathbb{N}
^{\ast }$. Then we have $x_{s}-x_{r}\in k^{r}%
\mathbb{Z}
^{n}$\ for $r<s$. If there exists no $x\in 
\mathbb{Z}
^{n}$ such that $x-x_{r}\in k^{r}%
\mathbb{Z}
^{n}$ for each $r$, then $\mathcal{E}$ satisfies (C). If $\mathcal{E}$
satisfies (C), then $\mathcal{E}$ is locally isomorphic to each $\widehat{%
\mathcal{A}}_{R}$ which satisfies (C).\bigskip

\noindent \textbf{Proof.} First suppose that there exist $r<s$ such that $%
x_{s}-x_{r}\notin k^{r}%
\mathbb{Z}
^{n}$. Then we have $x\sim _{\mathcal{E}}x+l^{r}y$ for $x\in 
\mathbb{Z}
^{n}-(x_{r}+k^{r}%
\mathbb{Z}
^{n})$ and\ $y\in 
\mathbb{Z}
^{n}$. We also have $x\sim _{\mathcal{E}}x+l^{s}y$ for $x\in 
\mathbb{Z}
^{n}-(x_{s}+k^{s}%
\mathbb{Z}
^{n})$ and\ $y\in 
\mathbb{Z}
^{n}$. It follows $x\sim _{\mathcal{E}}x+l^{s}y$ for each $x\in \mathcal{E}$
and each\ $y\in 
\mathbb{Z}
^{n}$, since $%
\mathbb{Z}
^{n}=(%
\mathbb{Z}
^{n}-(x_{r}+k^{r}%
\mathbb{Z}
^{n}))\cup (%
\mathbb{Z}
^{n}-(x_{s}+k^{s}%
\mathbb{Z}
^{n}))$.

By Proposition 1.2, there exist\ $t>s$,\ and $v,w\in k^{s}%
\mathbb{Z}
^{n}-k^{t}%
\mathbb{Z}
^{n}$\ such that $w-v\in l^{s}%
\mathbb{Z}
^{n}$\ and $v\nsim _{\mathcal{A}^{t}}w$. As $\Omega _{t,x_{t}}(\mathcal{E}%
)\cong \mathcal{A}^{t}$, there also exist $x,y\in x_{t}+(k^{s}%
\mathbb{Z}
^{n}-k^{t}%
\mathbb{Z}
^{n})$\ such that $y-x\in l^{s}%
\mathbb{Z}
^{n}$\ and $x\nsim _{\mathcal{E}}y$, whence a contradiction.

Now suppose that there exists no $x\in 
\mathbb{Z}
^{n}$ such that $x-x_{r}\in k^{r}%
\mathbb{Z}
^{n}$ for each $r$. Then, for each $x\in 
\mathbb{Z}
^{n}$ and each $h\in 
\mathbb{N}
^{\ast }$, there exists $s\in 
\mathbb{N}
^{\ast }$ such that $\mathrm{P}_{x,h}\cap (x_{s}+k^{s}%
\mathbb{Z}
^{n})=\emptyset $, since, otherwise, some $y\in \mathrm{P}_{x,h}$\ would
belong to infinitely many $x_{t}+k^{t}%
\mathbb{Z}
^{n}$, and therefore belong to all of them since $x_{t}\in x_{u}+k^{t}%
\mathbb{Z}
^{n}$ for $t<u$.\ For such an $s$, each isomorphism $\sigma :\Omega
_{s,x_{s}}(\mathcal{E})\rightarrow \mathcal{A}^{s}$ induces an isomorphism
from $\mathcal{E}\upharpoonright \mathrm{P}_{x,h}=\Omega _{s,x_{s}}(\mathcal{%
E})\upharpoonright \mathrm{P}_{x,h}$ to $\mathcal{A}^{s}\upharpoonright 
\mathrm{P}_{\sigma (x),h}=\widehat{\mathcal{A}}_{R}\upharpoonright \mathrm{P}%
_{\sigma (x),h}$.

It follows that $\mathcal{E}$ is locally isomorphic to each $\widehat{%
\mathcal{A}}_{R}$ which satisfies the local isomorphism property.
Consequently, $\mathcal{E}$ satisfies the local isomorphism property, and in
particular satisfies (C).

Now suppose that $\mathcal{E}$ satisfies (C) and that there exists $x\in 
\mathbb{Z}
^{n}$ such that $x-x_{r}\in k^{r}%
\mathbb{Z}
^{n}$ for each $r$. Then, for each $\widehat{\mathcal{A}}_{R}$ which
satisfies the local isomorphism property, in order to prove that $\mathcal{E}
$ is locally isomorphic to $\widehat{\mathcal{A}}_{R}$, il suffices to show
that, for each $h\in 
\mathbb{N}
^{\ast }$, there exists $z\in 
\mathbb{Z}
^{n}$ such that $\mathcal{E}\upharpoonright \mathrm{P}_{x,h}\cong \widehat{%
\mathcal{A}}_{R}\upharpoonright \mathrm{P}_{z,h}$.

We consider $y\in 
\mathbb{Z}
^{n}-\mathrm{P}_{x,h}$\ such that $\mathcal{E}\upharpoonright \mathrm{P}%
_{y,h}\cong \mathcal{E}\upharpoonright \mathrm{P}_{x,h}$\ and $s\in 
\mathbb{N}
^{\ast }$\ such that $\mathrm{P}_{y,h}\cap (x+k^{s}%
\mathbb{Z}
^{n})=\emptyset $. For each isomorphism $\sigma :\Omega _{s,x}(\mathcal{E}%
)=\Omega _{s,x_{s}}(\mathcal{E})\rightarrow \mathcal{A}^{s}$, we have

\noindent $\mathcal{E}\upharpoonright \mathrm{P}_{x,h}\cong \mathcal{E}%
\upharpoonright \mathrm{P}_{y,h}=\Omega _{s,x}(\mathcal{E})\upharpoonright 
\mathrm{P}_{y,h}\cong \mathcal{A}^{s}\upharpoonright \mathrm{P}_{\sigma
(y),h}=\widehat{\mathcal{A}}_{R}\upharpoonright \mathrm{P}_{\sigma (y),h}$.~~%
$\blacksquare $\bigskip

\noindent \textbf{Remark.} Suppose that there exists $x\in 
\mathbb{Z}
^{n}$ such that $x-x_{s}\in k^{s}%
\mathbb{Z}
^{n}$ for each $s\in 
\mathbb{N}
^{\ast }$. If there exists a sequence $(y_{s})_{s\in 
\mathbb{N}
^{\ast }}\subset 
\mathbb{Z}
^{n}-\left\{ x\right\} $ such that $y_{s}\in x+l^{s}%
\mathbb{Z}
^{n}$ and $y_{s}\sim _{\mathcal{E}}x$\ for each $s\in 
\mathbb{N}
^{\ast }$, then $\mathcal{E}$ satisfies (C) since $\mathcal{E}%
\upharpoonright \mathrm{P}_{x,k^{s}-1}\cong \mathcal{E}\upharpoonright 
\mathrm{P}_{y_{s},k^{s}-1}$\ for each $s\in 
\mathbb{N}
^{\ast }$.\ Otherwise, it is possible to change the connexions of the
segments of $\mathcal{E}$ with endpoint $x$ so that $(y_{s})_{s\in 
\mathbb{N}
^{\ast }}$ exists.\bigskip

\noindent \textbf{Remark.} If $\mathcal{E}$ is locally isomorphic to some $%
\widehat{\mathcal{A}}_{R}$ which satisfies (C), then $\mathcal{E}$ satisfies
the local isomorphism property like $\widehat{\mathcal{A}}_{R}$.\bigskip

The two examples below concern folding curves in $\mathbb{R}^{2}$. The
example in $\mathbb{R}^{3}$ that we give in Section 2 has some common points
with the first of them.

For each $s\in 
\mathbb{N}
^{\ast }$, an $s$-folding curve is obtained by folding $s$ times a strip of
paper in $2$, each time possibly left or right, then unfolding it with $\pi
/2$ angles.

Figure 1 shows a set $\mathcal{A}$ of equivalent curves which satisfies (H)
for $(2,4)$.\ For each $s\in 
\mathbb{N}
^{\ast }$, the curves of $\mathcal{A}^{s}$\ are called \emph{positive} $2s$%
-folding curves or $2s$-\emph{dragon} curves: in each of them, all the
foldings are done in the same direction.

Figure 2 shows the curves of $\mathcal{A}^{2}$ with endpoint $0$. We see
that $\widehat{\mathcal{A}}$ is the union of $4$ curves with endpoint $0$
which are equivalent modulo positive isometries. It follows from [6] that
each of the $2$ ways to realize the connexions in $0$ gives an $\widehat{%
\mathcal{A}}_{R}$ which satisfies the local isomorphism property. It
consists of $2$ curves without endpoints.

We define a fractal as follows: For each $s\in 
\mathbb{N}
^{\ast }$, we consider the $2s$-dragon curve $D_{s}\in \mathcal{A}^{s}$ with
endpoints $0$, $2^{s}e_{1}$; we denote by $F_{s}$ the image by a homothety
of center $0$ and ratio $1/2^{s}$ of the union of the Voronoi tiles
associated to the segments of $D_{s}$.\ The fractal is the limit of the sets 
$F_{s}$.

Figure 3 shows a set $\mathcal{B}$ of equivalent curves which also satisfies
(H) for $(2,4)$. For each $s\in 
\mathbb{N}
^{\ast }$, the curves of $\mathcal{B}^{s}$\ are \emph{alternate} $2s$%
-folding curves, obtained by folding $2s$ times a strip of paper in $2$,
alternatively left and right.

Figure 4 shows the curves of $\mathcal{B}^{2}$ which are contained in $\left[
-4,+4\right] ^{2}$, and in particular the $4$ curves which have $0$ as an
endpoint. We see that, contrary to the preceding case, $\widehat{\mathcal{B}}
$ is not the union of $4$ curves with endpoint $0$. It follows from [6]
that, also in this case, each of the $2$ ways to realize the connexions in $%
0 $ gives a $\widehat{\mathcal{B}}_{R}$ which satisfies the local
isomorphism property. It consists of $6$ curves without endpoints.

We show in [6] that each covering $\widehat{\mathcal{A}}_{R}$ or $\widehat{%
\mathcal{B}}_{R}$ is locally isomorphic to a covering which consists of one
curve. The same property will be true for the coverings considered
below.\bigskip

\textbf{2. Construction of the example in }$%
\mathbb{R}
^{3}$\textbf{.}\bigskip

From now on, we represent each element of $%
\mathbb{R}
^{3}$\ as a triple $(a,b,c)$. We write $0$ instead of $(0,0,0)$ when there
is no ambiguity. For the sake of brevity, we write $\overline{n}$\ instead
of $-n$\ for each $n\in 
\mathbb{N}
^{\ast }$.

For each curve $D$, we consider the curve $\overline{D}$ obtained by
reversing the order of the segments. We write $D_{1}=([0,e_{1}])$. For $%
s\geq 1$, we write $D_{s+1}=(D_{s},\left[ e_{s},e_{s}+e_{s+1}\right]
,e_{s+1}+\overline{D}_{s})$ and $C_{s+1}=(D_{s+1},\left[ e_{s+1},2e_{s+1}%
\right] )$.

We have $C_{2}=(\left[ 0,e_{1}\right] ,\left[ e_{1},e_{1}+e_{2}\right] ,%
\left[ e_{1}+e_{2},e_{2}\right] ,\left[ e_{2},2e_{2}\right] )$ and

\noindent $C_{3}=(\left[ 0,e_{1}\right] ,\left[ e_{1},e_{1}+e_{2}\right] ,%
\left[ e_{1}+e_{2},e_{2}\right] ,\left[ e_{2},e_{2}+e_{3}\right] ,$

$\ \ \ \left[ e_{2}+e_{3},e_{1}+e_{2}+e_{3}\right] ,\left[
e_{1}+e_{2}+e_{3},e_{1}+e_{3}\right] ,\left[ e_{1}+e_{3},e_{3}\right] ,\left[
e_{3},2e_{3}\right] )$.

The positive $2$-folding curve considered in the example above is equivalent
to $C_{2}$. Similarly, in order to construct our example in dimension $3$,
we are going to use the curve

\noindent $C=(\left[ 0,\overline{e}_{2}\right] ,\left[ \overline{e}_{2},%
\overline{e}_{2}+\overline{e}_{3}\right] ,\left[ \overline{e}_{2}+\overline{e%
}_{3},\overline{e}_{3}\right] ,\left[ \overline{e}_{3},e_{1}+\overline{e}_{3}%
\right] ,$

$\ \ \ \left[ e_{1}+\overline{e}_{3},e_{1}+\overline{e}_{2}+\overline{e}_{3}%
\right] ,\left[ e_{1}+\overline{e}_{2}+\overline{e}_{3},e_{1}+\overline{e}%
_{2}\right] ,\left[ e_{1}+\overline{e}_{2},e_{1}\right] ,\left[ e_{1},2e_{1}%
\right] )$,

\noindent which is equivalent to $C_{3}$ (see Figure 5). We can hope that
the curves $C_{s}$ for $s\geq 4$\ give analogous examples in higher
dimensions.

In each curve which is equivalent to $C$, we say that the image of $0$
(resp. $2e_{1}$) is an endpoint of \emph{type} $0$ (resp. $1$).

The curve $C$ belongs to a unique block which consists of $6$ disjoint
equivalent curves having $0$ as an endpoint of type $0$. The curves of this
block, are shown in Figure 6. They contain:

\noindent (i) the segments contained in $\mathrm{P}_{0,1}$, except the $12$
segments contained in its $6$ edges which have $(1,1,1)$ or $(\overline{1},%
\overline{1},\overline{1})$ as an endpoint;

\noindent (ii) the $6$ segments having one endpoint in the center of a face
of $\mathrm{P}_{0,1}$ and the other endoint outside $\mathrm{P}_{0,1}$%
.\bigskip

\noindent \textbf{Theorem 2.1.} The curve $C$ can be completed into $2$
equivalent coverings $\mathcal{C}_{1},\mathcal{C}_{2}$ of $\mathbb{R}^{3}$
by curves which are equivalent to $C$. Each $\mathcal{C}_{i}$ satisfies $4%
\mathbb{Z}
^{3}+\mathcal{C}_{i}=\mathcal{C}_{i}$. For each $\mathcal{C}_{i}$, each $%
x\in 
\mathbb{Z}
^{3}$ is an endpoint of type $0$ (resp. $1$) for $6$ curves if $x=(2k,2l,2m)$
with $k,l,m\in 
\mathbb{Z}
$ and $k+l+m$ even (resp. odd); no other $x\in 
\mathbb{Z}
^{3}$ is an endpoint.\bigskip

\noindent \textbf{Proof.} We denote by $\mathcal{B}_{0}$ the block mentioned
above. Figures 7 and 8 show the $2$ possible ways to place $3$ blocks $%
\mathcal{B}_{1}$, $\mathcal{B}_{2}$, $\mathcal{B}_{3}$, which are images of $%
\mathcal{B}_{0}$ by rotations with axes parallel to $e_{3}$ and angles $\pi
/2$, $\pi $, $3\pi /2$, in such a way that:

\noindent (i) the curves of $\mathcal{B}_{0}$, $\mathcal{B}_{1}$, $\mathcal{B%
}_{2}$, $\mathcal{B}_{3}$ are disjoint;

\noindent (ii) each of the $3$ edges of $\mathrm{P}_{0,1}$ having $(1,1,1)$
as an endpoint is covered by $2$ segments of $\mathcal{B}_{1}$, or $2$
segments of $\mathcal{B}_{2}$, or $2$ segments of $\mathcal{B}_{3}$.

\noindent In these figures, we only represent the edges of the cubes
associated to the $4$ blocks which are covered by segments of these blocks.

In both cases, the sets $4u+\mathcal{B}_{i}$ with $u\in 
\mathbb{Z}
^{3}$ and $0\leq i\leq 3$\ are disjoint. They form a covering of $\mathbb{R}%
^{3}$ since $\mathcal{B}_{0}\cup \mathcal{B}_{1}\cup \mathcal{B}_{2}\cup 
\mathcal{B}_{3}$ contains exactly $192$ segments like $[0,+4[^{3}$.

The covering associated to Figure 8 is the image of the covering associated
to Figure 7 by the positive isometry

\noindent $((0,0,0),(1,0,0),(0,1,0),(0,0,1))\rightarrow ((0,0,0),(0,%
\overline{1},0),(\overline{1},0,0),(0,0,\overline{1}))$.~~$\blacksquare $%
\bigskip

From now on, we consider the covering $\mathcal{C}$ associated to Figure 7.
By Theorem 2.1, $\mathcal{C}$ satisfies (H) for $\left( 2,4\right) $. We
define the sets $\mathcal{C}^{s}$ according to Section 1.

The block $\mathcal{B}_{0}$\ consists of the $6$ curves

\noindent $((0,0,0),(0,\overline{1},0),(0,\overline{1},\overline{1}),(0,0,%
\overline{1}),(1,0,\overline{1}),(1,\overline{1},\overline{1}),(1,\overline{1%
},0),(1,0,0),(2,0,0))$,

\noindent $((0,0,0),(0,0,1),(0,1,1),(0,1,0),(\overline{1},1,0),(\overline{1}%
,1,1),(\overline{1},0,1),(\overline{1},0,0),(\overline{2},0,0))$,

\noindent $((0,0,0),(0,0,\overline{1}),(\overline{1},0,\overline{1}),(%
\overline{1},0,0),(\overline{1},1,0),(\overline{1},1,\overline{1}),(0,1,%
\overline{1}),(0,1,0),(0,2,0))$,

\noindent $((0,0,0),(1,0,0),(1,0,1),(0,0,1),(0,\overline{1},1),(1,\overline{1%
},1),(1,\overline{1},0),(0,\overline{1},0),(0,\overline{2},0))$,

\noindent $((0,0,0),(\overline{1},0,0),(\overline{1},\overline{1},0),(0,%
\overline{1},0),(0,\overline{1},1),(\overline{1},\overline{1},1),(\overline{1%
},0,1),(0,0,1),(0,0,2))$,

\noindent $((0,0,0),(0,1,0),(1,1,0),(1,0,0),(1,0,\overline{1}),(1,1,%
\overline{1}),(0,1,\overline{1}),(0,0,\overline{1}),(0,0,\overline{2}))$.

The block $\mathcal{B}_{1}$\ consists of the $6$ curves

\noindent $%
((2,2,0),(2,3,0),(2,3,1),(2,2,1),(3,2,1),(3,3,1),(3,3,0),(3,2,0),(4,2,0))$,

\noindent $((2,2,0),(2,2,\overline{1}),(2,1,\overline{1}%
),(2,1,0),(1,1,0),(1,1,\overline{1}),(1,2,\overline{1}),(1,2,0),(0,2,0))$,

\noindent $((2,2,0),(3,2,0),(3,2,\overline{1}),(2,2,\overline{1}),(2,3,%
\overline{1}),(3,3,\overline{1}),(3,3,0),(2,3,0),(2,4,0))$,

\noindent $%
((2,2,0),(2,2,1),(1,2,1),(1,2,0),(1,1,0),(1,1,1),(2,1,1),(2,1,0),(2,0,0))$,

\noindent $%
((2,2,0),(2,1,0),(3,1,0),(3,2,0),(3,2,1),(3,1,1),(2,1,1),(2,2,1),(2,2,2))$,

\noindent $((2,2,0),(1,2,0),(1,3,0),(2,3,0),(2,3,\overline{1}),(1,3,%
\overline{1}),(1,2,\overline{1}),(2,2,\overline{1}),(2,2,\overline{2}))$.

The block $\mathcal{B}_{2}$\ consists of the $6$ curves

\noindent $((2,0,2),(2,0,3),(2,\overline{1},3),(2,\overline{1},2),(3,%
\overline{1},2),(3,\overline{1},3),(3,0,3),(3,0,2),(4,0,2))$,

\noindent $%
((2,0,2),(2,1,2),(2,1,1),(2,0,1),(1,0,1),(1,1,1),(1,1,2),(1,0,2),(0,0,2))$,

\noindent $%
((2,0,2),(1,0,2),(1,0,3),(2,0,3),(2,1,3),(1,1,3),(1,1,2),(2,1,2),(2,2,2))$,

\noindent $((2,0,2),(2,0,1),(3,0,1),(3,0,2),(3,\overline{1},2),(3,\overline{1%
},1),(2,\overline{1},1),(2,\overline{1},2),(2,\overline{2},2))$,

\noindent $%
((2,0,2),(3,0,2),(3,1,2),(2,1,2),(2,1,3),(3,1,3),(3,0,3),(2,0,3),(2,0,4))$,

\noindent $((2,0,2),(2,\overline{1},2),(1,\overline{1},2),(1,0,2),(1,0,1),(1,%
\overline{1},1),(2,\overline{1},1),(2,0,1),(2,0,0))$.

The block $\mathcal{B}_{3}$\ consists of the $6$ curves

\noindent $%
((0,2,2),(0,2,1),(0,3,1),(0,3,2),(1,3,2),(1,3,1),(1,2,1),(1,2,2),(2,2,2))$,

\noindent $((0,2,2),(0,1,2),(0,1,3),(0,2,3),(\overline{1},2,3),(\overline{1}%
,1,3),(\overline{1},1,2),(\overline{1},2,2),(\overline{2},2,2))$,

\noindent $%
((0,2,2),(0,2,3),(1,2,3),(1,2,2),(1,3,2),(1,3,3),(0,3,3),(0,3,2),(0,4,2))$,

\noindent $((0,2,2),(\overline{1},2,2),(\overline{1},2,1),(0,2,1),(0,1,1),(%
\overline{1},1,1),(\overline{1},1,2),(0,1,2),(0,0,2))$,

\noindent $((0,2,2),(0,3,2),(\overline{1},3,2),(\overline{1},2,2),(\overline{%
1},2,3),(\overline{1},3,3),(0,3,3),(0,2,3),(0,2,4))$,

\noindent $%
((0,2,2),(1,2,2),(1,1,2),(0,1,2),(0,1,1),(1,1,1),(1,2,1),(0,2,1),(0,2,0))$.

We write $\Delta (x)=2x$ for each $x\in 
\mathbb{Z}
^{3}$ and $\Delta (\mathcal{D})=\mathcal{C}\ast \mathcal{D}$ for each set $%
\mathcal{D}$ of disjoint self-avoiding curves. We denote by $\Gamma $ the
inverse operation: In each $\Delta (\mathcal{D})$, we first replace each
piece of curve with endpoints in $2%
\mathbb{Z}
^{3}$ consisting of $8$ segments by the segment which joins its endpoints,
then we apply a homothety of center $0$ and ratio $1/2$.

The operations $\Gamma $ and $\Delta \ $are analogous to the inflation and
the deflation considered for tilings. We have $\mathcal{C}^{n+1}=\Delta (%
\mathcal{C}^{n})$ for $n\in 
\mathbb{N}
$ and $\widehat{\mathcal{C}}=\Delta (\widehat{\mathcal{C}})$.

We note that, for $n\geq 2$, two consecutive segments of $\mathcal{C}^{n}$\
can be aligned or form a $\pi /2$ angle. Actually, one of the examples after
Theorem 2.7 below gives a curve of $\mathcal{C}^{2}$ such that $\left[ (2,2,%
\overline{1}),(2,2,0)\right] $ and $\left[ (2,2,0),(2,2,1)\right] $ are
consecutive segments.

For each $n\in 
\mathbb{N}
$, $\mathcal{C}^{n}$ is a covering of $\mathbb{R}^{3}$ by self-avoiding
curves which consist of $8^{n}$ segments. Each curve has endpoints $2^{n}u$, 
$2^{n}(u+e_{i})$ with $u\in 
\mathbb{Z}
^{3}$ and $i\in \left\{ 1,2,3\right\} $. We have $2^{n+1}%
\mathbb{Z}
^{3}+\mathcal{C}^{n}=\mathcal{C}^{n}$.

For any integers $m,n\geq 1$, $\mathcal{C}^{m+n}$ is obtained from $\mathcal{%
C}^{n}$ by applying a homothety of center $0$ and ratio $2^{m}$, then
remplacing each segment by the curve of $\mathcal{C}^{m}$ which has the same
endpoints.

We have $\mathcal{C}^{m+n}=\Delta ^{m}(\mathcal{C}^{n})$ and $\mathcal{C}%
^{n}=\Gamma ^{m}(\mathcal{C}^{m+n})$. The connexions between segments of $%
\mathcal{C}^{m+n}$\ in a point $x\in 
\mathbb{Z}
^{3}$ originate from $\mathcal{C}^{n}$\ via $\Delta ^{m}$\ if and only if $%
x\in 2^{m}%
\mathbb{Z}
^{3}$. By suppressing\ these connexions in $\mathcal{C}^{m+n}$, we get $%
\mathcal{C}^{m}$.\bigskip

\noindent \textbf{Remark.} The group consisting of isometries which
stabilize \textbf{$\mathcal{B}_{0}$} is generated by $\rho $ and $\sigma $
with $\rho (e_{1})=e_{2}$, $\rho (e_{2})=e_{3}$, $\rho (e_{3})=e_{1}$, $%
\sigma (e_{1})=\overline{e}_{2}$, $\sigma (e_{2})=\overline{e}_{1}$, $\sigma
(e_{3})=\overline{e}_{3}$.

The isometry $\rho $ stabilizes $\mathcal{C}$\ since it is linear and
permutes the sets \textbf{$\mathcal{B}_{i}$}. On the other hand, $\sigma $
does not stabilize $\mathcal{C}$. Actually, it sends the center $(2,0,2)$ of 
\textbf{$\mathcal{B}_{2}$} to $(0,\overline{2},\overline{2})$, which is the
center of $(0,\overline{4},\overline{4})+$\textbf{$\mathcal{B}_{3}$}, but it
sends $(1,\overline{1},3)$, which is not an endpoint of segments of curves
of \textbf{$\mathcal{B}_{2}$}, to $(1,\overline{1},\overline{3})$, which is
an endpoint of segments of curves of $(0,\overline{4},\overline{4})+$\textbf{%
$\mathcal{B}_{3}$}.

For each $s\in 
\mathbb{N}
^{\ast }$, the group consisting of isometries which stabilize $\mathcal{C}%
^{s}$ is generated by $\rho $ and the translations $x\rightarrow x+2^{s+1}u$
for $u\in \left\{ e_{1},e_{2},e_{3}\right\} $. The only isometries which
stabilize $\widehat{\mathcal{C}}$\ are $\mathrm{Id}$, $\rho $ and $\rho ^{2}$%
. For each pairing relation $R$ defined on $E$, $\rho $ and $\rho ^{2}$
stabilize $\widehat{\mathcal{C}}_{R}$ if and only if $R$ is invariant
through $\rho $.\bigskip

\noindent \textbf{Lemma 2.2.} Each curve of $\mathcal{C}^{4}$ contains a
curve of $\mathcal{C}^{2}$ which passes $3$ times through some point.\bigskip

\noindent \textbf{Proof.} First we show that the curve of $\mathcal{C}^{2}$
with endpoints $(0,0,0)$ and $(\overline{4},0,0)$ passes $3$ times through
some point. It is obtained from the curve

\noindent $((0,0,0),(0,0,1),(0,1,1),(0,1,0),(\overline{1},1,0),(\overline{1}%
,1,1),(\overline{1},0,1),(\overline{1},0,0),(\overline{2},0,0))$

\noindent of $\mathcal{C}$ by first applying a homothety of center $(0,0,0)$
and ratio $2$, then remplacing each of the segments

\noindent $\lbrack (0,0,0),(0,0,2)]$, $[(0,0,2),(0,2,2)]$, $%
[(0,2,2),(0,2,0)] $, $[(0,2,0),(\overline{2},2,0)]$,

\noindent $\lbrack (\overline{2},2,0),(\overline{2},2,2)]$, $[(\overline{2}%
,2,2),(\overline{2},0,2)]$, $[(\overline{2},0,2),(\overline{2},0,0)]$, $[(%
\overline{2},0,0),(\overline{4},0,0)]$

\noindent by the curve of $\mathcal{C}$ which has the same endpoints.

The segments $[(0,0,2),(0,2,2)]$,

\noindent $\lbrack (0,2,0),(\overline{2},2,0)]=(\overline{4}%
,0,0)+[(4,2,0),(2,2,0)]$,

\noindent $\lbrack (\overline{2},2,0),(\overline{2},2,2)]=(\overline{4}%
,0,0)+[(2,2,0),(2,2,2)]$

\noindent are respectively replaced by the curves

\noindent $((0,0,2),(0,1,2),(\overline{1},1,2),(\overline{1}%
,1,1),(0,1,1),(0,2,1),(\overline{1},2,1),(\overline{1},2,2),(0,2,2))$,

\noindent $((0,2,0),(\overline{1},2,0),(\overline{1},3,0),(\overline{1}%
,3,1),(\overline{1},2,1),(\overline{2},2,1),(\overline{2},3,1),(\overline{2}%
,3,0),(\overline{2},2,0))$,

\noindent $((\overline{2},2,0),(\overline{2},1,0),(\overline{1},1,0),(%
\overline{1},2,0),(\overline{1},2,1),(\overline{1},1,1),(\overline{2},1,1),(%
\overline{2},2,1),(\overline{2},2,2))$.

These curves are obtained by translations from the 4th curve of $\mathcal{B}%
_{3}$, the 1st curve of $\mathcal{B}_{1}$ and the 5th curve of $\mathcal{B}%
_{1}$; they all pass through the point $(\overline{1},2,1)$.

Now, for each $i\in \left\{ 1,2,3\right\} $ and each $x\in 8%
\mathbb{Z}
^{3}$, let us denote by $C_{i,x}$ the curve of $\mathcal{C}^{2}$ which joins 
$x-4e_{i}$ and $x$. We have just proved that $C_{1,0}$ passes $3$ times
through some point. The same property is true for $C_{2,0}$ and $C_{3,0}$
since they are images of $C_{1,0}$ by $\rho $ and $\rho ^{2}$. It is also
true for each curve $C_{i,x}$ since it is an image of $C_{i,0}$ by a
translation.

It suffices to show that each curve $A\in \mathcal{C}^{4}$ contains a curve $%
C_{i,x}$. For that purpose, we consider $2$ parallel segments $\left[
x,x+8e_{i}\right] $, $\left[ x+8e_{j},x+8e_{i}+8e_{j}\right] $, with $x\in 8%
\mathbb{Z}
^{3}$ and $i\neq j$, which have the same endpoints as $2$ curves of $%
\mathcal{C}^{3}$ contained in $A$. Then, the segments $\left[
x+4e_{i},x+8e_{i}\right] $, $\left[ x+4e_{i}+8e_{j},x+8e_{i}+8e_{j}\right] $
have the same endpoints as the curves $C_{i,x+8e_{i}}$ and $%
C_{i,x+8e_{i}+8e_{j}}$. One of these curves is necessarily contained in $A$%
.~~$\blacksquare $\bigskip

For each integer $s\geq 1$, we denote by $\mathcal{D}_{s}$ the set of all
curves of $\mathcal{C}^{s}$ which have $0$\ as an endpoint.\bigskip

\noindent \textbf{Proposition 2.3.} $\mathcal{D}_{2}$ covers $\mathrm{P}%
_{0,1}$.\bigskip

\noindent \textbf{Proof.} The curves of $\mathcal{D}_{1}$ contain the
segments contained in $\mathrm{P}_{0,1}$, except the $12$ segments contained
in the $6$ edges having $(1,1,1)$ or $(\overline{1},\overline{1},\overline{1}%
)$ as an endpoint.

Each of these $12$ segments is contained in a curve with endpoints $2x,2y$,
where $\left[ x,y\right] $ is a segment with $x\in \{(u,v,w)\in 
\mathbb{Z}
^{3}\mid \left| u\right| +\left| v\right| +\left| w\right| =1\}$ and $y\in
\{(u,v,w)\in 
\mathbb{Z}
^{3}\mid \left| u\right| +\left| v\right| +\left| w\right| =2\}$. This curve
is contained in a curve of $\mathcal{D}_{2}$ since its image $\left[ x,y%
\right] $ by $\Gamma $ belongs to a curve of $\mathcal{D}_{1}$.~~$%
\blacksquare $\bigskip

\noindent \textbf{Proposition 2.4.} For each integer $s\geq 4$, $\mathcal{D}%
_{s}$ covers $\mathrm{P}_{0,2^{s-3}}$. In particular, $\widehat{\mathcal{C}}$
is the union of $6$ curves with endpoint $0$.\bigskip

\noindent \textbf{Proof.} We write $\left\| (u,v,w)\right\| =\sup (\left|
u\right| ,\left| v\right| ,\left| w\right| )$ for $(u,v,w)\in 
\mathbb{R}
^{3}$ and $d(x,y)=\left\| y-x\right\| $ for $x,y\in 
\mathbb{R}
^{3}$. For each $x\in 
\mathbb{Z}
^{3}$ and each\ set $\mathcal{E}$ of curves, we write $d(x,\mathcal{E})=\inf
\{d(x,y)\mid y$ point of a curve of $\mathcal{E}\}$. We show by induction on 
$s$ that $d(0,\mathcal{C}^{s}-\mathcal{D}_{s})\geq 2^{s-3}+1$ for $s\geq 3$.

We have $d(0,\mathcal{C}^{3}-\mathcal{D}_{3})\geq 2$ since, by Proposition
2.3, for each $A\in \mathcal{C}^{3}-\mathcal{D}_{3}$, no segment of $\Gamma
(A)$ is contained in $\mathrm{P}_{0,1}$, and therefore no segment of $A$ is
contained in the interior of $\mathrm{P}_{0,2}$.

For $s\geq 3$, as $\mathcal{C}^{s+1}-\mathcal{D}_{s+1}=\Delta (\mathcal{C}%
^{s}-\mathcal{D}_{s})$, we have $d(0,\mathcal{C}^{s+1}-\mathcal{D}%
_{s+1})\geq 2d(0,\mathcal{C}^{s}-\mathcal{D}_{s})-1$. Consequently, $d(0,%
\mathcal{C}^{s}-\mathcal{D}_{s})\geq 2^{s-3}+1$ implies $d(0,\mathcal{C}%
^{s+1}-\mathcal{D}_{s+1})\geq 2^{s-2}+1$.~~$\blacksquare $\bigskip

\noindent \textbf{Proposition 2.5.} For each integer $s\geq 8$, each curve
of $\mathcal{C}^{s}$ covers some $\mathrm{P}_{x,2^{s-8}}$. In particular,
each curve of $\widehat{\mathcal{C}}$ covers arbitrarily large cubes.\bigskip

\noindent \textbf{Proof.} By Lemma 2.2, each curve of $\mathcal{C}^{4}$
passes $3$ times through some point. Consequently, each curve $A\in \mathcal{%
C}^{s}$ passes $3$ times through some point $x\in 2^{s-4}%
\mathbb{Z}
^{3}$.

The translation $0\rightarrow x$ sends $\mathcal{D}_{s-5}$ to the set of
curves of $\mathcal{C}_{{}}^{s-5}$ which have $x$ as an endpoint. The curves
of that set are contained in $A$. They cover $\mathrm{P}_{x,2^{s-8}}$ since $%
\mathcal{D}_{s-5}$ covers $\mathrm{P}_{0,2^{s-8}}$ by Proposition 2.4.~~$%
\blacksquare $\bigskip

Lemma 2.6 below will be used to prove Theorem 2.7.

For each set $\mathcal{E}$ of curves, each $x\in 
\mathbb{Z}
^{3}$ and each $s\in 
\mathbb{N}
^{\ast }$, we denote by $\Sigma _{s,x}(\mathcal{E})$\ the set of curves
obtained from $\mathcal{E}$ as follows: for each curve $D\in \mathcal{E}$\
with $x\in D$, we take the curve or the $2$ curves with endpoint $x$
contained in $D$ and having maximal length, with this length limited to $%
8^{s}$.

For each $x\in 2^{s+1}%
\mathbb{Z}
^{3}$, we have $\Sigma _{s,x}(\mathcal{C}^{s})=x+\Sigma _{s,0}(\mathcal{C}%
^{s})$ and $\Sigma _{s,x}(\mathcal{C}^{s})$ is the set of curves of $%
\mathcal{C}^{s}$ which have $x$ as an endpoint.\bigskip

\noindent \textbf{Lemme 2.6.} For each $s\in 
\mathbb{N}
^{\ast }$, there exists no $x\in 
\mathbb{Z}
^{3}-2^{s+1}%
\mathbb{Z}
^{3}$ such that $\Sigma _{s,x}(\mathcal{C}^{s})\prec x+\Sigma _{s,0}(%
\mathcal{C}^{s})$.\bigskip

\noindent \textbf{Proof.} We first observe that it suffices to show this
result for $s=1$. Actually, if we have $1\leq t<s$ for the largest integer $%
t $ such that $x\in 2^{t+1}%
\mathbb{Z}
^{3}$, then\ $\Sigma _{s,x}(\mathcal{C}^{s})\prec x+\Sigma _{s,0}(\mathcal{C}%
^{s})$ implies $\Sigma _{s-t,x/2^{t}}(\mathcal{C}^{s-t})\prec x/2^{t}+\Sigma
_{s-t,0}(\mathcal{C}^{s-t})$ since $\Gamma ^{t}(\mathcal{C}^{s})=\mathcal{C}%
^{s-t}$ and $\Gamma ^{t}(x)=x/2^{t}$; it follows $\Sigma _{1,x/2^{t}}(%
\mathcal{C})=\Sigma _{1,x/2^{t}}(\mathcal{C}^{s-t})\prec x/2^{t}+\Sigma
_{1,0}(\mathcal{C}^{s-t})=x/2^{t}+\Sigma _{1,0}(\mathcal{C})$ with $%
x/2^{t}\in 2%
\mathbb{Z}
^{3}-4%
\mathbb{Z}
^{3}$.

Now let us consider $x\in 
\mathbb{Z}
^{3}-4%
\mathbb{Z}
^{3}$ such that $\Sigma _{1,x}(\mathcal{C})\prec x+\Sigma _{1,0}(\mathcal{C}%
) $.

First suppose that $x\in 
\mathbb{Z}
^{3}-2%
\mathbb{Z}
^{3}$. Consider a curve $B\in \mathcal{C}$ which passes through $x$ and
write $B=B_{1}\cup B_{2}$ with $B_{1},B_{2}\in \Sigma _{1,x}(\mathcal{C})$.
Then we see that $B_{1}$ and $B_{2}$ cannot simultaneously be images of
subcurves of $C$ by positive isometries which send $0$ to $x$, which
contradicts $\Sigma _{1,x}(\mathcal{C})\prec x+\Sigma _{1,0}(\mathcal{C})$.

Now suppose that $x\in 2%
\mathbb{Z}
^{3}-4%
\mathbb{Z}
^{3}$.\ Then, there exists $y=\alpha _{1}e_{1}+\alpha _{2}e_{2}+\alpha
_{3}e_{3}$ with $\left( \alpha _{1},\alpha _{2},\alpha _{3}\right) \in
\left\{ 0,2\right\} ^{3}-\left\{ \left( 0,0,0\right) \right\} $\ such that $%
x\in y+4%
\mathbb{Z}
^{3}$.\ We have $\Sigma _{1,y}(\mathcal{C})\cong \Sigma _{1,x}(\mathcal{C})$
and therefore $\Sigma _{1,y}(\mathcal{C})\prec y+\Sigma _{1,0}(\mathcal{C})$.

If $\alpha _{1}+\alpha _{2}+\alpha _{3}=2$ or $\alpha _{1}+\alpha
_{2}+\alpha _{3}=6$, then\ $y$ is an endpoint of type $1$ for the curves of $%
\Sigma _{1,y}(\mathcal{C})$, which contradicts $\Sigma _{1,y}(\mathcal{C}%
)\prec y+\Sigma _{1,0}(\mathcal{C})$.\ If $\alpha _{1}+\alpha _{2}+\alpha
_{3}=4$, then $\Sigma _{1,y}(\mathcal{C})$ is the image of $\Sigma _{1,0}(%
\mathcal{C})$ by the composition of a translation and a rotation of angle $%
\pi /2$, $\pi $ or $3\pi /2$, which also contradicts $\Sigma _{1,y}(\mathcal{%
C})\prec y+\Sigma _{1,0}(\mathcal{C})$.~~$\blacksquare $\bigskip

In $\widehat{\mathcal{C}}$, the connexions between the segments which have
the same endpoint $x$ are defined for $x\neq 0$. Only the connexions in $0$
remain to be defined.

For $x\in 2%
\mathbb{Z}
^{3}$ and $u,v\in E$, we write $\theta _{x}(u)=v$ if the curve of $\mathcal{C%
}$ which connects $x$ and $x+2u$ contains the segment $[x,x+v]$.

For $x\in 
\mathbb{Z}
^{3}-\left\{ 0\right\} $ and $u,v\in E$, we write $R_{x}(u,v)$ if the
segments $\left[ x,x+u\right] $ and $\left[ x,x+v\right] $ are connected. As 
$\Gamma (\widehat{\mathcal{C}})=\widehat{\mathcal{C}}$, we have $R_{x}(u,v)$%
\ if and only if $R_{2x}(\theta _{2x}(u),\theta _{2x}(v))$.

As $\mathcal{C}$ is stabilized by the translations $u\rightarrow u+4v$ for $%
v\in 
\mathbb{Z}
^{3}$, we have $R_{x}=R_{x+4y}$ for $x\in 
\mathbb{Z}
^{3}-2%
\mathbb{Z}
^{3}$ and $y\in 
\mathbb{Z}
^{3}$, and $\theta _{x}=\theta _{x+4y}$ for $x\in 2%
\mathbb{Z}
^{3}$ and $y\in 
\mathbb{Z}
^{3}$. It follows $R_{x}=R_{x+8y}$ for $x\in 2%
\mathbb{Z}
^{3}-4%
\mathbb{Z}
^{3}$ and $y\in 
\mathbb{Z}
^{3}$.

We consider the permutation $\theta $ of $E$ with $\theta (e_{1})=\overline{e%
}_{2}$, $\theta (e_{2})=\overline{e}_{3}$, $\theta (e_{3})=\overline{e}_{1}$
and $\theta ^{2}=\mathrm{Id}$. We have $\theta _{x}=\theta $ for each $x\in 4%
\mathbb{Z}
^{3}$.

For each relation $R$ defined on $E$ and any $u,v\in E$, we write $R^{\theta
}(u,v)$ if and only if $R(\theta (u),\theta (v))$. For each $x\in 2%
\mathbb{Z}
^{3}$, as $\theta _{2x}=\theta $, we have $R_{x}=R$ (resp. $%
R_{x}=R_{{}}^{\theta }$) if and only if $R_{2x}=R_{{}}^{\theta }$ (resp. $%
R_{2x}=R$).

For each pairing relation $R$ defined on $E$, we consider the set of
complete curves $\widehat{\mathcal{C}}_{R}$\ obtained from $\widehat{%
\mathcal{C}}$\ by defining the connexions in $0$ according to $R$. We say
that $R$ satisfies (P) if there exists $x\in \left\{ \overline{2}%
,0,2,4\right\} ^{3}-\left\{ 0,4\right\} ^{3}$ such that $R_{x}=R$ or $%
R_{x}=R_{{}}^{\theta }$.\bigskip

\noindent \textbf{Theorem 2.7.} $\widehat{\mathcal{C}}_{R}$ satisfies the
local isomorphism property if and only if $R$ satisfait (P).\bigskip

\noindent \textbf{Proof.} If $\widehat{\mathcal{C}}_{R}$\ satisfies the
local isomorphism property, then there exists $u\in 
\mathbb{Z}
^{3}-P_{0,2}$ such that $\widehat{\mathcal{C}}_{R}\upharpoonright \mathrm{P}%
_{0,2}\cong \widehat{\mathcal{C}}_{R}\upharpoonright \mathrm{P}_{u,2}$. \ It
follows $R_{u}=R$ and $\Sigma _{1,u}(\mathcal{C})\prec \Sigma _{1,0}(%
\mathcal{C})$. By Lemma 2.6, we have $u\in 2^{s+1}%
\mathbb{Z}
^{3}-2^{s+2}%
\mathbb{Z}
^{3}$ for the largest integer $s$ such that $\Sigma _{s,u}(\mathcal{C}%
^{s})\prec \Sigma _{s,0}(\mathcal{C}^{s})$.

There exists $v\in 2%
\mathbb{Z}
^{3}-4%
\mathbb{Z}
^{3}$\ such that $u=2^{s}v$. We have $R_{v}=R$ (resp. $R_{v}=R^{\theta }$)
if $s$ is even (resp. odd). The same property is true if we remplace $v$ by
any $w\in 
\mathbb{Z}
^{3}$\ such that $w-v\in 8%
\mathbb{Z}
^{3}$. There exists such a $w$ which belongs to $\left\{ \overline{2}%
,0,2,4\right\} ^{3}$.

Conversely, let us suppose that there exists $x\in 2%
\mathbb{Z}
^{3}-4%
\mathbb{Z}
^{3}$\ such that $R_{x}=R$\ (resp.\ $R_{x}=R^{\theta }$). Then we have $%
R_{2^{s}x+2^{s+3}y}=R$ for each\ even (resp. odd) integer $s$ and each $y\in 
\mathbb{Z}
^{3}$. We also have $\Sigma _{s-1,2^{s}x+2^{s+3}y}(\mathcal{C}^{s-1})\cong
\Sigma _{s-1,0}(\mathcal{C}^{s-1})$ for any such $s,y$. Now, it follows from
Proposition 2.4 that $\widehat{\mathcal{C}}_{R}$\ satisfies the local
isomorphism property.~~$\blacksquare $\bigskip

\noindent \textbf{Examples.} Here, we give examples of values of $R_{x}$\
for $x\in \left\{ \overline{2},0,2,4\right\} ^{3}-\left\{ 0,4\right\} ^{3}$.
By Theorem 2.7, they are also values of $R$ such that $\widehat{\mathcal{C}}%
_{R}$\ satisfies the local isomorphism property.

For $x=(a,b,c)$ with $a+b+c$ odd, as $2x$ is an endpoint of type $1$, we
have $\theta _{2x}=\mathrm{Id}$ and therefore $R_{2x}=R_{x}$. In particular,

\noindent $R_{1,1,1}=\left\langle (e_{1},\overline{e}_{3}),(e_{2},\overline{e%
}_{1}),(e_{3},\overline{e}_{2})\right\rangle $ and

\noindent $R_{\overline{1},\overline{1},\overline{1}}=R_{1,0,0}=R_{\overline{%
1},0,0}=R_{0,1,0}=R_{0,\overline{1},0}=R_{0,0,1}=R_{0,0,\overline{1}}$

$\ \ \ \ \ \ =\left\langle (e_{1},\overline{e}_{2}),(e_{2},\overline{e}%
_{3}),(e_{3},\overline{e}_{1})\right\rangle $

\noindent respectively imply $R_{2,2,2}=\left\langle (e_{1},\overline{e}%
_{3}),(e_{2},\overline{e}_{1}),(e_{3},\overline{e}_{2})\right\rangle $ and

\noindent $R_{\overline{2},\overline{2},\overline{2}}=R_{2,0,0}=R_{\overline{%
2},0,0}=R_{0,2,0}=R_{0,\overline{2},0}=R_{0,0,2}=R_{0,0,\overline{2}}$

$\ \ \ \ \ \ =\left\langle (e_{1},\overline{e}_{2}),(e_{2},\overline{e}%
_{3}),(e_{3},\overline{e}_{1})\right\rangle $.

We also have $R_{1,1,0}=\left\langle (e_{1},\overline{e}_{3}),(e_{2},e_{3}),(%
\overline{e}_{1},\overline{e}_{2})\right\rangle $ and

\noindent $\theta _{2,2,0}:\left( e_{1},\overline{e}_{1},e_{2},\overline{e}%
_{2},e_{3},\overline{e}_{3}\right) \rightarrow \left( e_{2},\overline{e}%
_{3},e_{1},e_{3},\overline{e}_{2},\overline{e}_{1}\right) $.

\noindent It follows $R_{2,2,0}=\left\langle (e_{1},\overline{e}_{2}),(e_{2},%
\overline{e}_{1}),(e_{3},\overline{e}_{3})\right\rangle $. We note that the
aligned segments $\left[ (2,2,\overline{1}),(2,2,0)\right] $ and $\left[
(2,2,0),(2,2,1)\right] $\ are connected. By applying the isometry $\rho $,
we obtain $R_{0,2,2}=\left\langle (e_{1},\overline{e}_{1}),(e_{2},\overline{e%
}_{3}),(e_{3},\overline{e}_{2})\right\rangle $ and $R_{2,0,2}=\left\langle
(e_{1},\overline{e}_{3}),(e_{2},\overline{e}_{2}),(e_{3},\overline{e}%
_{1})\right\rangle $.\bigskip

The following theorem gives a characterization of the sets of curves which
are locally isomorphic to the sets $\widehat{\mathcal{C}}_{R}$ for $R$
satisfying (P).

By Lemma 2.6, for $s,t\in 
\mathbb{N}
^{\ast }$, $s\leq t$ and $x,y\in 
\mathbb{Z}
^{3}$, we have $x+\mathcal{C}^{s}\prec y+\mathcal{C}^{t}$\ if and only if $%
y-x\in 2^{s+1}%
\mathbb{Z}
^{3}$. For each sequence $X=(x_{s})_{s\in 
\mathbb{N}
^{\ast }}\subset 
\mathbb{Z}
^{3}$ with $x_{s+1}-x_{s}\in 2^{s+1}%
\mathbb{Z}
^{3}$ for $s\in 
\mathbb{N}
^{\ast }$, we denote by $\widehat{\mathcal{C}}_{X}$ the inductive limit of
the sets $x_{s}+\mathcal{C}^{s}$.\bigskip

\noindent \textbf{Theorem 2.8.} For each relation $R$ which satisfies (P),
any set of curves $\mathcal{E}$ is locally isomorphic to $\widehat{\mathcal{C%
}}_{R}$ if and only if it satisfies 1) or 2) below:

\noindent 1) $\mathcal{E}\cong \widehat{\mathcal{C}}_{S}$ for a relation $S$
which satisfies (P);

\noindent 2) $\mathcal{E}\cong \widehat{\mathcal{C}}_{X}$ for a sequence $%
X=(x_{s})_{s\in 
\mathbb{N}
^{\ast }}$\ such that $\bigcap_{s\in 
\mathbb{N}
^{\ast }}x_{s}+2^{s+1}%
\mathbb{Z}
^{3}=\emptyset $.\bigskip

\noindent \textbf{Proof.} It follows from Theorems 1.4 and 2.7 that $%
\mathcal{E}$ is locally isomorphic to $\widehat{\mathcal{C}}_{R}$ if 1) is
true, and also if 2) is true since $\bigcap_{s\in 
\mathbb{N}
^{\ast }}x_{s}+2^{s+1}%
\mathbb{Z}
^{3}=\bigcap_{s\in 
\mathbb{N}
^{\ast }}x_{s}+2^{s}%
\mathbb{Z}
^{3}$.

By Theorem 1.3, if $\mathcal{E}$ is locally isomorphic to $\widehat{\mathcal{%
C}}_{R}$, then, for each $s\in 
\mathbb{N}
^{\ast }$, there exists $x_{s}\in 
\mathbb{Z}
^{3}$ such that $(\Omega _{s,x_{s}}(\mathcal{E}),x_{s})\cong (\mathcal{C}%
^{s},0)$. For each $s\in 
\mathbb{N}
^{\ast }$, $(\Omega _{s,x_{s}}(\mathcal{E}),x_{s})\cong (\Omega _{s,x_{s+1}}(%
\mathcal{E}),x_{s+1})\cong (\mathcal{C}^{s},0)$ implies $x_{s+1}-x_{s}\in
2^{s+1}%
\mathbb{Z}
^{3}$ by Lemma 2.6.

The property 2) is true if $\bigcap_{s\in 
\mathbb{N}
^{\ast }}x_{s}+2^{s+1}%
\mathbb{Z}
^{3}$ is empty. If $\bigcap_{s\in 
\mathbb{N}
^{\ast }}x_{s}+2^{s+1}%
\mathbb{Z}
^{3}$ contains an element $x$, then the inductive limit of the sets $\Omega
_{s,x}(\mathcal{E})=\Omega _{s,x_{s}}(\mathcal{E})$ is the image of $%
\widehat{\mathcal{C}}$ by the translation $u\rightarrow x+u$. It follows
that $\mathcal{E}$ is the image of $\widehat{\mathcal{C}}_{S}$ by this
translation for a relation $S$; this relation satisfies (P) by Theorem 2.7.~~%
$\blacksquare $\bigskip

\noindent \textbf{Theorem 2.9.} There exists a covering of $\mathbb{R}^{3}$
by a self-avoiding curve which is locally isomorphic to the sets $\widehat{%
\mathcal{C}}_{R}$ for $R$ satisfying (P).\bigskip

\noindent \textbf{Proof.} We fix $R$ which satisfies (P). We consider the
cubes $\mathrm{P}_{x,h}$.

We first show that, for each cube $K$, there exist an integer $s$, a curve $%
A\in \mathcal{C}^{s}$ and a cube $L$ covered by $A$ such that $%
(L,A\upharpoonright L)\cong (K,\widehat{\mathcal{C}}_{R}\upharpoonright K)$.
As $\widehat{\mathcal{C}}_{R}$\ satisfies the local isomorphism property,
there exists $h\in 
\mathbb{N}
^{\ast }$ such that each $\mathrm{P}_{x,h}$ contains a cube $M$ with $(M,%
\widehat{\mathcal{C}}_{R}\upharpoonright M)\cong (K,\widehat{\mathcal{C}}%
_{R}\upharpoonright K)$. By Proposition 2.5, for such an $h$, there exists $%
s\in 
\mathbb{N}
^{\ast }$ such that each $A\in \mathcal{C}^{s}$ covers some $\mathrm{P}%
_{x,h} $, and therefore covers a cube $L$ with $(L,A\upharpoonright L)\cong
(K,\widehat{\mathcal{C}}_{R}\upharpoonright K)$.

As each bounded curve is contained in a cube, it follows that there exists a
sequence $(K_{s},A_{s})_{s\in 
\mathbb{N}
}$ such that, for each $s\in 
\mathbb{N}
$, $K_{s}$\ is a cube and $A_{s}$ is a\ curve covering $K_{s}$, isomorphic
to a curve of some $\mathcal{C}^{m}$, contained in the interior of $K_{s+1}$
and restriction of $A_{s+1}$ to the union of the segments of $A_{s}$.\ The
union of the curves $A_{s}$ is a self-avoiding curve covering $\mathbb{R}%
^{3} $ and locally isomorphic to $\widehat{\mathcal{C}}_{R}$ by Theorem
2.8.~~$\blacksquare $\bigskip

\noindent \textbf{Proposition 2.10.} Let $\mathcal{D}$ be a finite set of
curves such that $\Gamma ^{n}(\mathcal{D})$\ is defined for each $n\in 
\mathbb{N}
$. Then there exist $n\in 
\mathbb{N}
$ and $x\in 
\mathbb{Z}
^{3}$ such that each curve of $\Gamma ^{n}(\mathcal{D})$\ contains a point
of $\mathrm{P}_{x,1}$.\bigskip

\noindent \textbf{Proof.} We define the distance $d$ as in the proof of
Proposition 2.5.

For each $z\in \mathbb{R}^{3}$ and each $D\in \mathcal{D}$, we have

\noindent $d(\Gamma ^{n+1}(z),\Gamma ^{n+1}(D))\leq (1/2)(d(\Gamma
^{n}(z),\Gamma ^{n}(D))+1)$

\noindent for each $n\in 
\mathbb{N}
$, and therefore $d(\Gamma ^{n}(z),\Gamma ^{n}(D))<3/2$ for $n$ large
enough. Consequently, there exist $y\in \mathbb{R}^{3}$ and $n\in 
\mathbb{N}
$ such that $d(y,\Gamma ^{n}(D))<3/2$ for each $D\in \mathcal{D}$.

Now, let us consider $x\in 
\mathbb{Z}
^{3}$ such that $d(x,y)\leq 1/2$. Then, for each $D\in \mathcal{D}$, we have 
$d(x,\Gamma ^{n}(D))\leq d(x,y)+d(y,\Gamma ^{n}(D))<3/2+1/2=2$, and
therefore $d(x,\Gamma ^{n}(D))\leq 1$\ since $d(x,\Gamma ^{n}(D))$ is an
integer.~~$\blacksquare $\bigskip

\noindent \textbf{Corollary 2.11.} Let $\mathcal{D}$ be a set of curves
without endpoints such that $\Gamma ^{n}(\mathcal{D})$\ is defined for each $%
n\in 
\mathbb{N}
$. Then $\mathcal{D}$ consists of at most $27$ curves.\bigskip

\noindent \textbf{Proof.} It suffices to show this result for $\mathcal{D}$
finite. By Proposition 2.10, there exist $n\in 
\mathbb{N}
$ and $x\in 
\mathbb{Z}
^{3}$ such that each curve of $\Gamma ^{n}(\mathcal{D})$\ contains a point
of $%
\mathbb{Z}
^{3}\cap \mathrm{P}_{x,1}$. Any such curve contains at least $2$ segments
with one endpoint in $(%
\mathbb{Z}
^{3}\cap \mathrm{P}_{x,1})-\left\{ x\right\} $ and the other endpoint
outside $\mathrm{P}_{x,1}$.

The set $(%
\mathbb{Z}
^{3}\cap \mathrm{P}_{x,1})-\left\{ x\right\} $ consists of:

\noindent - the $8$ vertices of $\mathrm{P}_{x,1}$, each of them being an
endpoint of $3$ segments whose other endpoint is outside $\mathrm{P}_{x,1}$;

\noindent - the middles of the $12$ edges of $\mathrm{P}_{x,1}$, each of
them being an endpoint of $2$ segments whose other endpoint is outside $%
\mathrm{P}_{x,1}$;\ \ 

\noindent - the middles of the $6$ faces of $\mathrm{P}_{x,1}$, each of them
being an endpoint of $1$ segment whose other endpoint is outside $\mathrm{P}%
_{x,1}$.

As the curves of $\Gamma ^{n}(\mathcal{D})$\ are disjoint, their number is
at most $(1/2)(3.8+2.12+6)=27$.~~$\blacksquare $\bigskip

\noindent \textbf{Remark.} The bound given in Corollary 2.11 is not a priori
optimal, contrary to those which were given in [6] and [8] for sets of
folding curves in $%
\mathbb{R}
^{2}$.\bigskip

The following result is analogous to [6, Th. 3.10] and [8, Th. 2.1], which
concern folding curves in $%
\mathbb{R}
^{2}$.\bigskip

\noindent \textbf{Proposition 2.12.} Let $\mathcal{E},\mathcal{F}$\ be sets
of curves which are locally isomorphic to the sets $\widehat{\mathcal{C}}%
_{R} $ for $R$ satisfying (P), and let $A$ be an unbounded curve contained
in a curve of $\mathcal{E}$ and in a curve of $\mathcal{F}$. If there exist
an $x\in 
\mathbb{Z}
^{3}$ and a relation $S$ satisfying (P) such that $(\mathcal{E},x)\cong (%
\widehat{\mathcal{C}}_{S},0)$, then $\mathcal{E}$ and $\mathcal{F}$ can only
differ by the connexions between their segments which have $x$ as an
endpoint.\ Otherwise, we necessarily have $\mathcal{E}=\mathcal{F}$. \bigskip

\noindent \textbf{Proof.} For each $x\in 
\mathbb{Z}
^{3}$, we denote by $\Omega _{x,\mathcal{E}}$\ (resp. $\Omega _{x,\mathcal{F}%
}$) the set of curves with endpoint $x$ and length $8$ contained in $%
\mathcal{E}$ (resp. $\mathcal{F}$). We write $M_{x}=%
\{x+2k_{1}e_{1}+2k_{2}e_{2}+2k_{3}e_{3}\mid k_{1},k_{2},k_{3}\in 
\mathbb{Z}
$ and $k_{1}+k_{2}+k_{3}$\ even$\}$.

We consider $x\in 
\mathbb{Z}
^{3}$\ such that, for each $y\in M_{x}$, the curves of $\Omega _{y,\mathcal{E%
}}$\ are equivalent to $C$ and $y$ is an endpoint of type $0$ for these
curves. Then, because of the existence of $A$, the same property is true for 
$\mathcal{F}$. Moreover, if $y\in M_{x}$ is an endpoint of a segment of $A$,
then we have $\Omega _{y,\mathcal{F}}=\Omega _{y,\mathcal{E}}$ since $\Omega
_{y,\mathcal{E}}$ and $\Omega _{y,\mathcal{F}}$ are blocks and contain the
same curve which is contained in $A$. It follows that $\Omega _{y,\mathcal{F}%
}=\Omega _{y,\mathcal{E}}$ for each $y\in M_{x}$.

Now, let us consider the sets of curves $\Gamma (\mathcal{E})$ and $\Gamma (%
\mathcal{F})$\ obtained from $\mathcal{E}$ and $\mathcal{F}$ by remplacing
each curve of $\cup _{y\in M_{x}}\Omega _{y,\mathcal{E}}=\cup _{y\in
M_{x}}\Omega _{y,\mathcal{F}}$ by the segment which has the same endpoints,
and applying a homothety of ratio $1/2$. Let us denote by $\Gamma (A)$ the
largest curve $B$ such that $\Delta (B)\subset A$.

Then, $\Gamma (\mathcal{E})$ and $\Gamma (\mathcal{F})$\ are locally
isomorphic to the sets $\widehat{\mathcal{C}}_{R}$ for $R$ satisfying (P).
Moreover, $\Gamma (A)$ is an unbounded curve contained in a curve of $\Gamma
(\mathcal{E})$ and in a curve of $\Gamma (\mathcal{F})$.

Now, the proposition follows from Theorem 2.8.~~$\blacksquare $\bigskip

Similar to the case of dragon curves in $\mathbb{R}^{2}$,\ the present
construction can be used to define a fractal. We consider the curves $%
C_{s}=\Delta ^{s}(C)$\ for $s\in 
\mathbb{N}
^{\ast }$. We denote by $F$\ the limit of the sets $F_{s}=(1/2^{s})(\cup
_{S\in C_{s}}V(S))$ where each $V(S)$\ is the Voronoi tile associated to $S$%
. Using Proposition 2.5, we see that $F$ is the closure of its interior.

We note that, similar to the fractal constructed from dragon curves, $F$
satisfies a property of \emph{autosimilarity}: For each $s\in 
\mathbb{N}
^{\ast }$, there exist $8$ disjoint curves $C_{s,1},\ldots ,C_{s,8}$ which
are images of $C_{s}$ by positive isometries and such that $%
C_{s+1}=C_{s,1}\cup \cdots \cup C_{s,8}$. Let us take them in such a way
that they are consecutive, starting from $0$. Then $F$ is the union of the $%
8 $ nonoverlapping sets $F_{i}$ obtained as limits of the sets $%
F_{i,s}=(1/2^{s})(\cup _{S\in C_{i,s}}V(S))$. Each $F_{i}$ is the image of $%
F $ by a positive similarity of ratio $1/2$.\bigskip

\begin{center}
\textbf{References\bigskip }
\end{center}

\noindent \lbrack 1] S.I. Ben-Abraham, A. Quandt and D. Shapira,
Multidimensional paperfolding systems, Acta Cryst. A 69 (2013), 123-130.

\noindent \lbrack 2] M. Dekking, Paperfolding morphisms, paperfolding
curves, and fractal tiles, Theoret. Comput. Sci. 414 (2012), 20-37.

\noindent \lbrack 3] C. Goodman-Strauss, C.S. Kaplan, J.S. Myers, and D.
Smith, A chiral aperiodic monotile, Combinatorial theory 4 (1) (2024).

\noindent \lbrack 4] M. Mend\`{e}s France and O. Shallit, Wire Bending, J.
Combinatorial Theory 50 (1989), 1-23.

\noindent \lbrack 5] F. Oger, Algebraic and model-theoretic properties of
tilings,Theoret. Comput. Sci. 319 (2004), 103-126.

\noindent \lbrack 6] F. Oger, Paperfolding sequences, paperfolding curves
and local isomorphism, Hiroshima Math. Journal 42 (2012), 37-75.

\noindent \lbrack 7] F. Oger, The number of paperfolding curves in a
covering of the plane, Hiroshima Math. Journal 47 (2017), 1-14.

\noindent \lbrack 8] F. Oger, Coverings of the plane by self-avoiding curves
which satisfy the local isomorphism property, arXiv:2310.19364.\bigskip

\begin{center}
\includegraphics[scale=0.52]{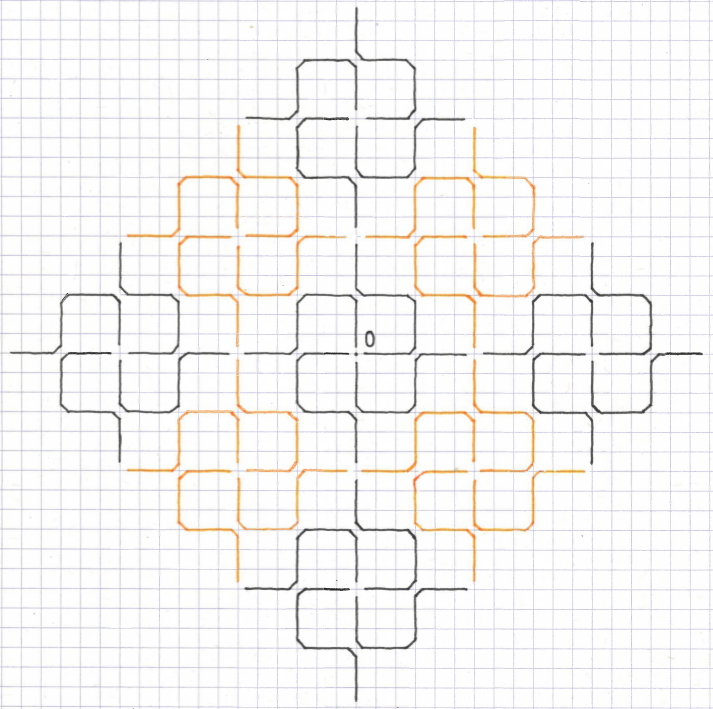}

\medskip Figure 1
\end{center}

\bigskip

\begin{center}
\includegraphics[scale=0.52]{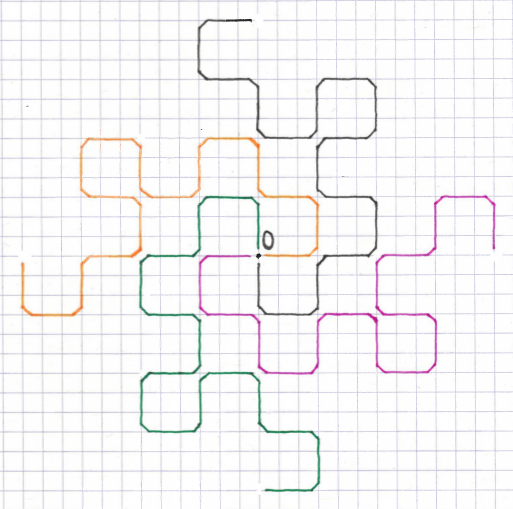}

\medskip Figure 2

\bigskip

\includegraphics[scale=0.52]{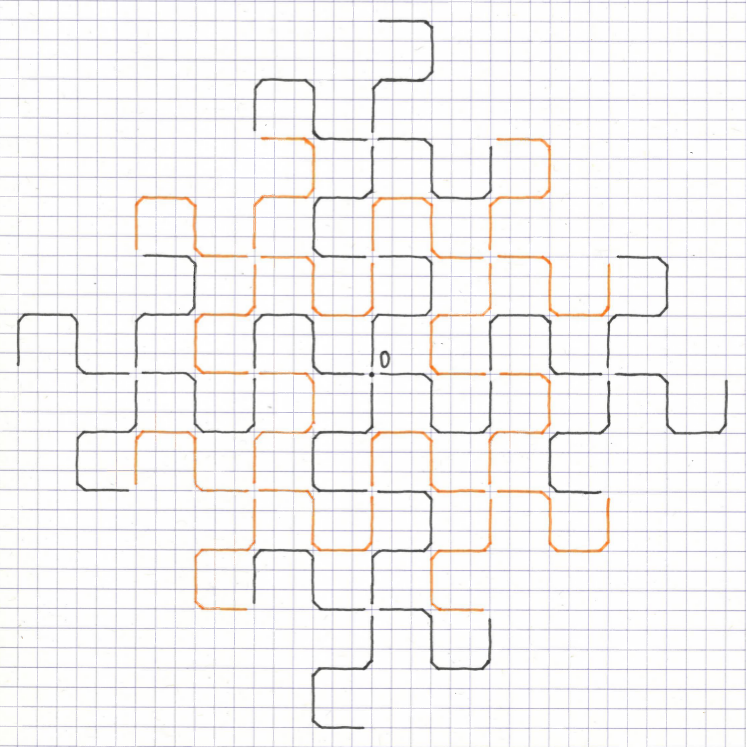}

\medskip Figure 3
\end{center}

\bigskip

\begin{center}
\includegraphics[scale=0.52]{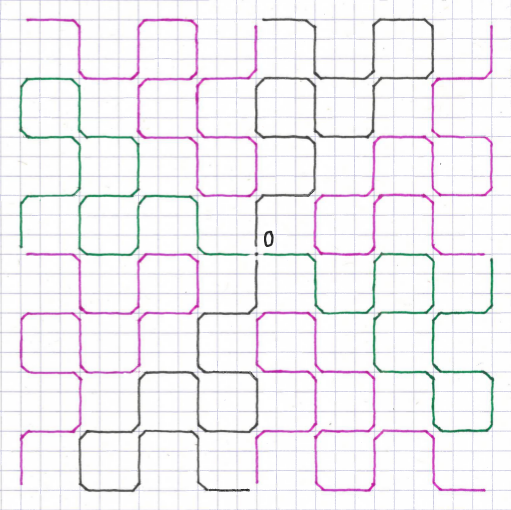}

\medskip Figure 4
\end{center}

\bigskip

\begin{center}
\includegraphics[scale=0.52]{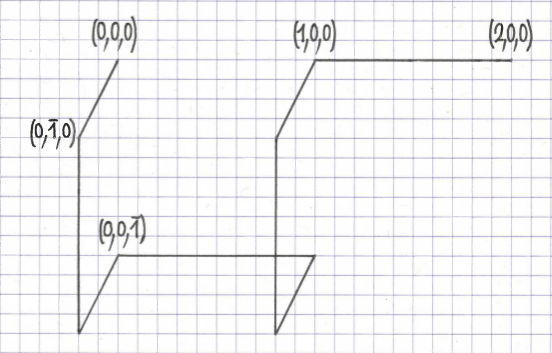}

\medskip Figure 5
\end{center}

\bigskip

\begin{center}
\includegraphics[scale=0.52]{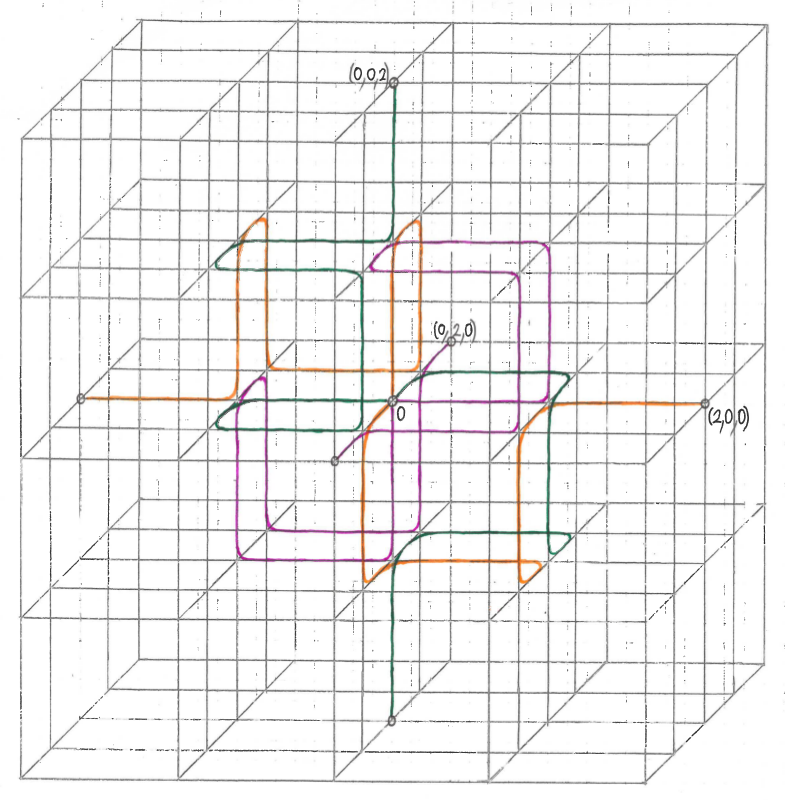}

\medskip Figure 6
\end{center}

\bigskip

\begin{center}
\includegraphics[scale=0.40]{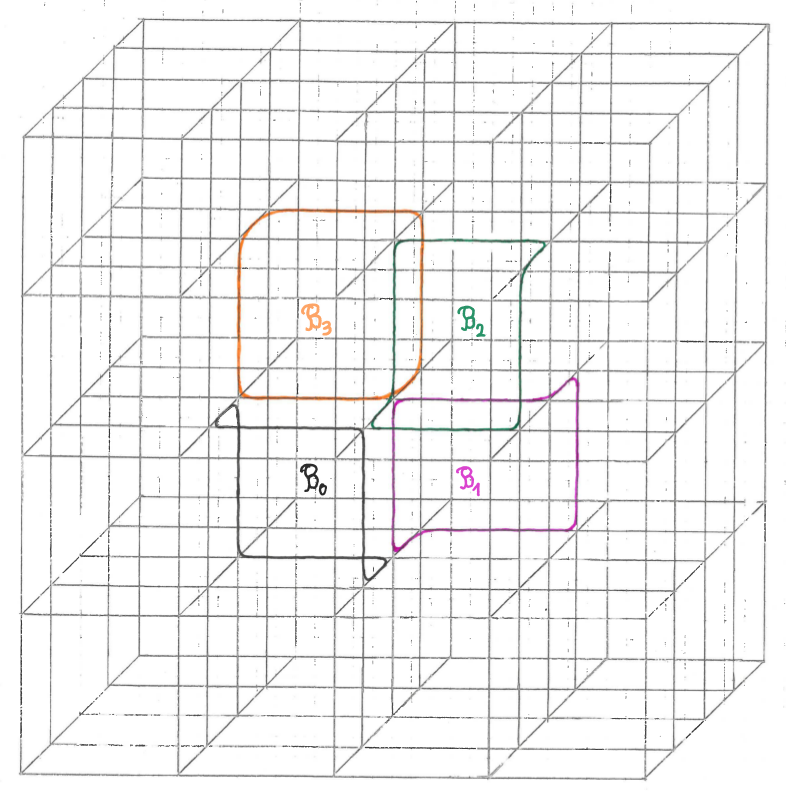}

\medskip Figure 7
\end{center}

\bigskip

\begin{center}
\includegraphics[scale=0.40]{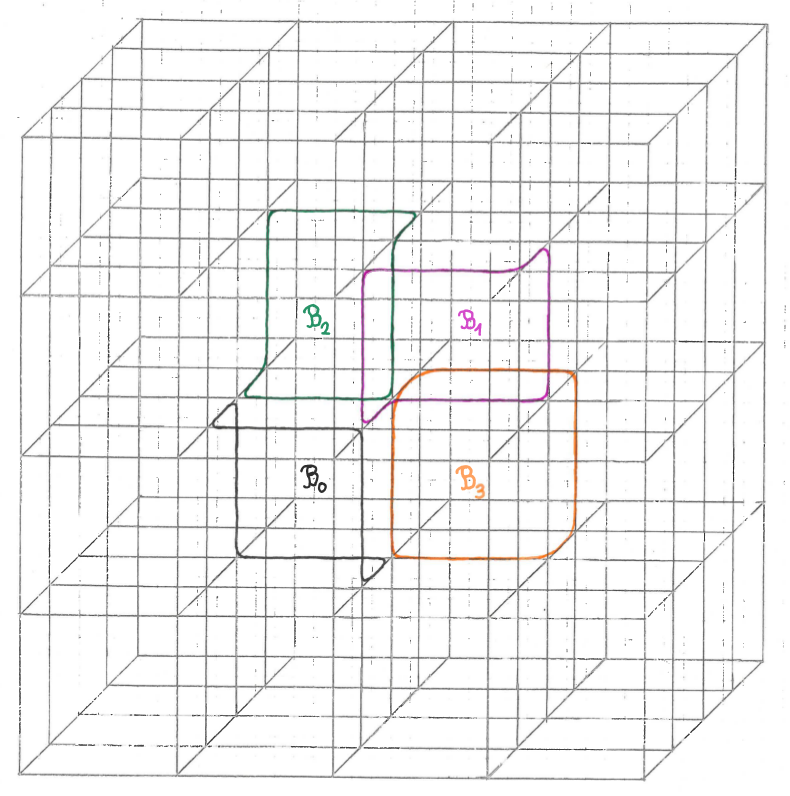}

\medskip Figure 8\bigskip 
\end{center}

Francis OGER

UFR de Math\'{e}matiques, Universit\'{e} Paris Cit\'{e}

B\^{a}timent Sophie Germain

8 place Aur\'{e}lie Nemours

75013 Paris

France

E-mail: oger@math.univ-paris-diderot.fr

\vfill

\end{document}